\documentclass[11pt]{amsart}

\usepackage[T1]{fontenc}
\usepackage{lmodern}
\usepackage[utf8]{inputenc}
\usepackage{amsmath,amssymb,amsthm}
\usepackage[hidelinks]{hyperref}
\usepackage{mathrsfs}
\usepackage{textcase}
\usepackage{tikz}
\usepackage{enumitem}
\usepackage{geometry}
\usepackage[backend=biber,sorting=nyt,doi=true,url=true]{biblatex}
\AtEveryBibitem{%
  \iffieldundef{doi}{}{%
    \clearfield{url}%
    \clearfield{urldate}%
    \clearfield{eprint}%
  }%
}

\newtheorem{theorem}{Theorem}[section]
\newtheorem*{unnumberedtheorem}{Theorem}
\newtheorem{problem}[theorem]{Problem}
\newtheorem{lemma}[theorem]{Lemma}
\newtheorem*{unnumberedlemma}{Lemma}
\newtheorem{proposition}[theorem]{Proposition}
\newtheorem{corollary}[theorem]{Corollary}

\theoremstyle{definition}
\newtheorem{definition}[theorem]{Definition}
\theoremstyle{remark}

\DeclareMathOperator{\Exp}{exp}

\numberwithin{equation}{section}

\title{Pseudocompact MAD families under
  \texorpdfstring{\NoCaseChange{$\mathfrak b=\mathfrak c$}}{b=c}}
\author{Vinicius de Oliveira Rodrigues}
\author{Fabrício Yohei Ynone}
\address{Instituto de Matemática, Estatística e Ciência da Computação,
  Universidade de São Paulo, Rua do Matão, 1010, 05508-090 São Paulo, SP, Brazil}
\email[V. de Oliveira Rodrigues]{vinior@ime.usp.br}
\email[F. Y. Ynone]{fabricio.ynone@usp.br}
\date{}

\begin{document}

\begin{abstract}
We prove that if the bounding number $\mathfrak b$ equals the cardinality
of the continuum $\mathfrak c$, then every infinite almost disjoint family
of size less than $\mathfrak c$ can be extended to a maximal almost disjoint
family whose Isbell--Mrówka space has a pseudocompact Vietoris hyperspace.
\end{abstract}

\maketitle

\section{Introduction}

The Vietoris hyperspace of a topological space is a space whose
points are the nonempty closed subsets of the original space.
Its topology was introduced by Leopold Vietoris in 1922 \cite{Vietoris1922}.
For compact metric spaces, it coincides with the topology induced by
the Hausdorff metric.

A central theme in hyperspace theory is
to understand how topological properties of a space are reflected
in its hyperspace and vice versa.
A classical result due to Vietoris states that
the Vietoris hyperspace of a topological space is compact if
and only if the original space is compact
\cite{Vietoris1922}.
Thus, it is natural to ask what are the relations between
generalizations of compactness between a space and its hyperspace.

In this direction, Ginsburg \cite{Ginsburg1975} investigated
countable compactness and pseudocompactness of Vietoris hyperspaces.
He proved that, for a Hausdorff space $X$:
\begin{enumerate}[label=(\roman*)]
    \item \label{item:countably-compact-filter} If $X$ is $p$-compact for some free ultrafilter $p \in \omega^*$, then $\Exp(X)$ is $p$-compact.
    \item \label{item:countably-compact-finite} If $\Exp(X)$ is countably compact, then $X^n$ is countably compact for every $n \in \omega$.
    \item \label{item:countably-compact-powers} If all the powers of $X$ are countably compact, then $\Exp(X)$ is countably compact.
\end{enumerate}
He also proved the following for a Tychonoff space $X$.
\begin{enumerate}[label=(\roman*')]
    \item \label{item:pseudocompact-filter} If $X$ is $p$-pseudocompact for some $p \in \omega^*$, then $\Exp(X)$ is $p$-pseudocompact.
    \item \label{item:pseudocompact-finite} If $\Exp(X)$ is pseudocompact, then $X^n$ is pseudocompact for every $n \in \omega$.
\end{enumerate}

Notice that \ref{item:pseudocompact-filter} and
\ref{item:pseudocompact-finite} are the pseudocompact analogues
of \ref{item:countably-compact-filter} and
\ref{item:countably-compact-finite}, respectively.
Moreover, \ref{item:countably-compact-powers} follows from
\ref{item:countably-compact-filter}: if all powers of $X$ are
countably compact, then $X$ is $p$-compact for some
$p\in\omega^*$ by the Ginsburg--Saks theorem
\cite[Theorem~2.6]{GinsburgSaks1975}. Applying
\ref{item:countably-compact-filter} to this $p$
shows that $\Exp(X)$ is $p$-compact, and hence
countably compact.
Despite that, as shown in \cite{HrusakEtAl2007}, there exists a
Tychonoff space $X$ such that all powers of $X$ are pseudocompact,
while its hyperspace is not pseudocompact---namely, a subspace of
$\beta\omega$.

Recall that if the $\omega$-th power of a Tychonoff space $X$ is pseudocompact, then all the powers of $X$ are pseudocompact.
Ginsburg asked how the pseudocompactness of $X^\omega$ is related to
that of $\Exp(X)$ and which spaces $X$ have pseudocompact
hyperspaces \cite[Remark~3.2, p.~1399]{Ginsburg1975}.
In this paper, we further explore these questions in the context of Isbell--Mr\'owka spaces.

Given an infinite almost disjoint family $\mathcal A$ on $\omega$, one
associates to it a topological space $\Psi(\mathcal A)$ known as
the Isbell--Mr\'owka space of $\mathcal A$
(\cite{Mrowka1954} and \cite[5I]{GillmanJerison1960}).
The pseudocompactness of $\Psi(\mathcal A)$ is equivalent to the maximality of $\mathcal A$, and these spaces are Tychonoff.
So the class of Isbell--Mr\'owka spaces provides a natural setting to study the Ginsburg's questions, which becomes:

\begin{problem}
    What is the relation between the maximality of an infinite almost disjoint family $\mathcal A$ and the pseudocompactness of the Vietoris hyperspace $\Exp(\Psi(\mathcal A))$?
\end{problem}

In this paper, a \emph{pseudocompact MAD family} is an infinite maximal almost disjoint family $\mathcal A$ such that $\Exp(\Psi(\mathcal A))$ is pseudocompact.

Several previous results have been obtained.
The following theorem shows that a full equivalence between the maximality of an almost disjoint family and the pseudocompactness of the hyperspaces is equivalent to a known combinatorial cardinal.

\begin{proposition}[{\cite[Theorem~2.4]{GuzmanEtAl2022}}]
The following are equivalent:
    \begin{enumerate}[label=(\roman*)]
        \item Every MAD family is pseudocompact.
        \item $n(\omega^*)>\mathfrak c$.
    \end{enumerate}
\end{proposition}

Here, $n(\omega^*)$ is the Novák number of $\omega^*$, which is the least cardinal $\kappa$ such that $\omega^*$ can be covered by $\kappa$ nowhere dense sets.
It is independent of ZFC whether $n(\omega^*)>\mathfrak c$, e.g., it holds under $\mathfrak p=\mathfrak c$, but fails under $\mathfrak h<\mathfrak c$.
In particular, if $\mathfrak p=\mathfrak c$, then every MAD family is pseudocompact, and if $\mathfrak h<\mathfrak c$, then there exists a MAD family that is not pseudocompact --- two results that had previously appeared in \cite{HrusakEtAl2007}.

Thus, the existence of a non-pseudocompact MAD family is independent of ZFC.
On the other hand, the existence of a pseudocompact MAD family in ZFC alone is an open problem.

\begin{problem}
    Does there exist a pseudocompact MAD family in ZFC?
\end{problem}

This question appears in
\cite[Problem~4.3]{HrusakEtAl2007},
\cite[Problem~9]{Hrusak2014}, \cite[Question~8.9]{HernandezHrusak2018}, \cite[Question~1.3]{GuzmanEtAl2022},
\cite[Question~1.4]{CorralRodriguesFin2024},
and \cite[Question~2]{CorralRodrigues2024}.
The first author has also presented talks on partial results regarding this problem at the 56th Spring Topology
and Dynamical Systems Conference and at the 38th Summer Conference on
Topology and its Applications
\cite{RodriguesSTDC2023,RodriguesSUMTOPO2024}.

Besides $\mathfrak p=\mathfrak c$, previous work shows that $\diamondsuit(\mathfrak b)$ and $\mathfrak{ap}=\mathfrak c$ also imply the existence of a pseudocompact MAD family \cite{CorralRodrigues2024}.
Brendle introduced and studied the cardinal $\mathfrak{ap}$, which satisfies $\mathfrak p\leq\mathfrak{ap}\leq\mathfrak b$ \cite{Brendle1999}.
He also shows that the strict inequality $\mathfrak{ap}<\mathfrak b$ is consistent with ZFC.
We strengthen the previous result by proving generic existence under $\mathfrak b=\mathfrak c$ in the sense of \cite[Section~3]{GuzmanEtAl2022}.

\begin{theorem}\label{thm:main}
    If $\mathfrak b=\mathfrak c$, then every infinite almost disjoint family of size less than $\mathfrak c$ can be extended to a pseudocompact MAD family.
\end{theorem}
Together with the previous results, we conclude that in a very broad class of models of ZFC, there exists a pseudocompact MAD family.
However, the question of whether there exists a pseudocompact MAD family in ZFC remains open.

\section{Preliminaries}

For subsets of $\omega$, write $x\subseteq^*y$ if
$x\setminus y$ is finite, and $x=^*y$ if their symmetric
difference is finite. An \emph{almost disjoint} (AD) family
is a collection of infinite subsets of $\omega$ whose distinct
members have finite intersection.
A \emph{MAD family (maximal almost disjoint family)} is
an infinite AD family maximal under inclusion.
Write
$r\perp\mathcal A$ if $r\cap a$ is finite for every
$a\in\mathcal A$.

For $f,g\in\omega^\omega$, write $f\leq^*g$ if $f(n)\leq g(n)$ for all but finitely many $n$.
The bounding number $\mathfrak b$ is the least cardinality of a subset of $\omega^\omega$ with no common $\leq^*$-upper bound.
The dominating number $\mathfrak d$ is the least cardinality of a set $D\subseteq\omega^\omega$ such that for every $f\in\omega^\omega$ there exists $g\in D$ such that $f\leq^*g$.
The almost disjointness number $\mathfrak a$ is the least cardinality of a MAD family.
For background on these cardinal characteristics of the continuum, see \cite{Blass2010}.

Given an almost disjoint family $\mathcal A$, the Isbell--Mr\'owka
space associated to $\mathcal A$ is denoted by $\Psi(\mathcal A)$
and has underlying set $\omega\mathbin{\dot\cup}\mathcal A$.
It has the topology for which the natural numbers
are isolated, and a neighborhood base at $a\in\mathcal A$
consists of the sets $\{a\}\cup(a\setminus s)$,
$s\in[\omega]^{<\omega}$.
This is the strongest topology on $\Psi(\mathcal A)$
for which every point $a\in\mathcal A\subseteq\Psi(\mathcal A)$
is the limit of the increasing enumeration of $a$.
It is well known that $\Psi(\mathcal A)$ is Tychonoff,
locally compact, zero-dimensional and separable.
Notice also that every subset of $\mathcal A$ is closed in
$\Psi(\mathcal A)$.
For background on almost disjoint families and Isbell--Mr\'owka spaces, see
\cite{Hrusak2014} and \cite{HernandezHrusak2018}.

A space is pseudocompact if every continuous real-valued
function on it is bounded.

\begin{proposition}[{\cite{HrusakEtAl2007}}]
  The following are equivalent for an infinite almost disjoint family $\mathcal A$:
  \begin{enumerate}[label=(\roman*)]
  \item $\mathcal A$ is a MAD family.
  \item The space $\Psi(\mathcal A)$ is pseudocompact.
  \item The space $\Psi(\mathcal A)^\omega$ is pseudocompact.
  \end{enumerate}
\end{proposition}

For a topological space $X$,
$\Exp(X)$ denotes its nonempty closed subsets with the
Vietoris topology. A basic open set of $\Exp(X)$ is of the form
\[
 \langle U_0,\ldots,U_{m-1}\rangle
 =\left\{H\in\Exp(X):
    H\subseteq\bigcup_{i<m}U_i
    \text{ and }H\cap U_i\neq\varnothing\text{ for all }i<m\right\},
\]
where $0<m<\omega$ and each $U_i$ is open in $X$.
If $U\subseteq X$, we put
\[
  U^+=\langle U\rangle = \{H\in\Exp(X):H\subseteq U\}\quad \text{and}\quad
  U^-=\langle X, U\rangle = \{H\in\Exp(X):H\cap U\neq\varnothing\}.
\]

If $X$ is a regular space, then $\Exp(X)$ is Hausdorff \cite{Michael1951}.
Thus, for every almost disjoint family $\mathcal A$, the hyperspace
$\Exp(\Psi(\mathcal A))$ is Hausdorff.

We say that an AD family $\mathcal A$ is \emph{pseudocompact} if
$\Exp(\Psi(\mathcal A))$ is pseudocompact.

Following \cite{CorralRodrigues2024}, a \emph{fin-sequence}
is a sequence $C=(C_n)_{n\in\omega}$ of pairwise disjoint
nonempty finite subsets of $\omega$. For $x,X\subseteq\omega$
and a finite family $F$ of subsets of $\omega$, set
\begin{align*}
 \pi_C(x)&=\{n\in\omega:C_n\cap x\neq\varnothing\},\\
 E_C(X)&=\{n\in\omega:C_n\subseteq X\},\\
 J_C(F)&=\bigcap_{b\in F}\pi_C(b),\qquad J_C(\varnothing)=\omega.
\end{align*}
For a family $B$, write $\pi_C[B]=\{\pi_C(b):b\in B\}$.
A family of subsets of $\omega$ is \emph{centered} if every
finite intersection is infinite. In particular, the empty
family is centered.

We use the following criterion to study pseudocompact MAD families.

\begin{proposition}[Rodrigues--Tomita {\cite[Proposition~2.1]{RodriguesTomita2019}}]
\label{prop:finite-criterion}
For an AD family $\mathcal A$, the hyperspace
$\Exp(\Psi(\mathcal A))$ is pseudocompact if and only if
every fin-sequence has an accumulation point in it.
\end{proposition}

To produce accumulation points, we use the following accumulation criterion which appears implicitly in \cite[Lemma~2.2]{GuzmanEtAl2022}.
We include a proof here for the convenience of the reader.

\begin{lemma}\label{lem:vietoris}
Let $\mathcal A$ be an almost disjoint family, let $C$ be a fin-sequence and let 
$\mathcal B\subseteq\mathcal A$ be nonempty.
\begin{enumerate}[label=(\roman*)]
  \item\label{item:vietoris-points} Every accumulation point of the sequence $C:\omega\rightarrow\Exp(\Psi(\mathcal A))$
is a nonempty subset of $\mathcal A$.
  \item\label{item:vietoris-criterion} The set $\mathcal B$ is an accumulation point of $C$ if and only if for every $F\in[\mathcal B]^{<\omega}$ and every $X \subseteq \omega$ such that $\forall b\in \mathcal B\, b\subseteq^* X$, we have
  \[
    |J_C(F)\cap E_C(X)|=\omega.
  \]
\end{enumerate}

\end{lemma}

\begin{proof}
Suppose that $H$ is an accumulation point of $C$.
Assume by contradiction that $H$ is not a subset of $\mathcal A$ and fix $k \in H\cap \omega$.
As $C$ is a pairwise disjoint sequence,
\[
  |\{i \in \omega: C_i\in (\{k\})^{-}\}|\leq 1,
\]
a contradiction as $(\{k\})^{-}$ is a neighborhood of $H$ in $\Exp(\Psi(\mathcal A))$.

Now assume $\mathcal B$ is an accumulation point and fix $X$ and $F$
as in the statement, so the set $U=\mathcal B\cup X$ is open in
$\Psi(\mathcal A)$.
For each $b\in F$, let $V_b=\{b\}\cup b$, which is also open.
Then put
\[
 \mathcal U=U^+\cap\bigcap_{b\in F}V_b^-.
\]

Notice that $\{n \in \omega: C_n \in \mathcal U\} = J_C(F)\cap E_C(X)$, thus,
the latter must be infinite.

Conversely, let
$\langle U_0,\ldots,U_{m}\rangle$ be a basic neighborhood of $\mathcal B$.
Choose $b_i\in \mathcal B\cap U_i$ for each $i\leq m$, and put
$F=\{b_i:i\leq m\}$. Since $U_i$ is a neighborhood of $b_i$,
the set $s_i=b_i\setminus U_i$ is finite. Set
\[
 X=\left(\omega\cap\bigcup_{i\leq m}U_i\right)
       \setminus\bigcup_{i\leq m}s_i.
\]
Every $b\in \mathcal B$ belongs to some $U_i$, so $b\subseteq^*X$. Hence
$|J_C(F)\cap E_C(X)|=\omega$ by hypothesis.
Now it is straightforward to check that
\[
  J_C(F)\cap E_C(X)\subseteq \{n\in\omega:C_n\in\langle U_0,\ldots,U_m\rangle\},
\]
which shows that $\mathcal B$ is an accumulation point of $C$.
\end{proof}

An ideal $\mathcal I$ on $\omega$ is \emph{free} if $[\omega]^{<\omega}\subseteq \mathcal I$.
For an ideal $\mathcal I$, let $\mathcal I^+=\mathcal P(\omega)\setminus\mathcal I$.
A generator of a free ideal $\mathcal I$ is a set $\mathcal G\subseteq\mathcal I$ such that $\mathcal I=\{X\subseteq\omega:\exists G\in\mathcal G\ X\subseteq^*G\}$. We write $\chi(\mathcal I)$ for the least cardinality of such a generator.

The ideal generated by an AD family $\mathcal A$ and its complement are denoted by
\[
  \mathcal I(\mathcal A)
  =\left\{Y\subseteq\omega:
    (\exists K\in[\mathcal A]^{<\omega})\
      Y\subseteq^*\bigcup K\right\},
  \qquad
  \mathcal I^+(\mathcal A)=\mathcal P(\omega)\setminus\mathcal I(\mathcal A).
\]

Thus, if $\mathcal A$ is an infinite almost disjoint family, then $\chi(\mathcal I(\mathcal A))=|\mathcal A|$.
If $\mathcal A$ is a finite almost disjoint family, then $\chi(\mathcal I(\mathcal A))\leq \omega$.

The following proposition is stated in \cite[Remark~3.5]{BlassHrusakVerner2013},
and a proof is implicit.
For the convenience of the reader, we include a proof here.

\begin{proposition}
\label{prop:small-ideal}
Every free ideal $\mathcal I$ with
$\chi(\mathcal I)<\mathfrak d$ is a $P^+$-ideal:
every decreasing sequence $(T_n)_{n\in\omega}$ in
$\mathcal I^+$ has an $I\in\mathcal I^+$ such that
$I\subseteq^*T_n$ for all $n$.
\end{proposition}

\begin{proof}
Fix a generator $\mathcal B$ for $\mathcal I$ of cardinality
$\chi(\mathcal I)$ and $(T_n)_{n\in\omega}$ a decreasing sequence in $\mathcal I^+$.
For each $B\in\mathcal B$ and $n \in \omega$, we have $T_n\not\subseteq^*B$, so $T_n\setminus B$ is infinite.
Define $f_B\in\omega^\omega$ by
\[
 f_B(n)=\min\bigl((T_n\setminus B)\setminus n\bigr).
\]
Then $|\{f_B:B\in\mathcal B\}|\leq\chi(\mathcal I)<\mathfrak d$, so $\{f_B:B\in\mathcal B\}$ is not dominating.
Hence there exists
$g\in\omega^\omega$ such that $g\not\leq^*f_B$ for every
$B\in\mathcal B$. Put
\[
 I=\bigcup_{n\in\omega}\bigl(T_n\cap[n,g(n))\bigr).
\]
The sequence $(T_n)_{n\in\omega}$
is decreasing, and therefore, for each $k\in\omega$,
\[
 I\setminus T_k\subseteq
 \bigcup_{n<k}[n,g(n)),
\]
which is finite. Thus $I\subseteq^*T_k$ for every $k$.
Finally, fix $B\in\mathcal B$. There exist infinitely many $n$
such that $f_B(n)<g(n)$, and for each such $n$,
\[
 f_B(n)\in (T_n\setminus B)\cap[n,g(n))\subseteq I\setminus B.
\]
Since $f_B(n)\geq n$ for each $n$, we conclude that $I\setminus B$ is infinite.
Thus, for no $B\in\mathcal B$ is $I\subseteq^*B$, and therefore $I\in \mathcal I^+$.
\end{proof}

As a corollary, we obtain the following.

\begin{corollary}\label{cor:positive-pseudo}
Let $\mathcal H\subseteq[\omega]^\omega$ be centered,
$|\mathcal H|<\mathfrak d$, and let $(H_n)_{n\in\omega}$
be a sequence in $\mathcal H$. There exists $I\in[\omega]^\omega$
such that $I\subseteq^*H_n$ for all $n\in\omega$, and
\[
 |I\cap\bigcap\mathcal K|=\omega
 \quad\text{for every }\mathcal K\in[\mathcal H]^{<\omega}.
\]
\end{corollary}
\begin{proof}
Let $\mathcal F$ be the filter generated by $\mathcal H$ together
with the cofinite subsets of $\omega$. Since $\mathcal H$ is centered,
$\mathcal F$ is a proper free filter.
Let $\mathcal I=\mathcal F^*=\{\omega\setminus F: F \in \mathcal F\}$ be the dual ideal of $\mathcal F$.
Then $\mathcal I$ has a set of at most $|\mathcal H|+\omega<\mathfrak d$ generators so by Proposition~\ref{prop:small-ideal}, $\mathcal I$ is a $P^+$-ideal.
Notice that $\mathcal F\subseteq \mathcal I^+$.

As $\mathcal F$ contains all finite intersections of members of $\mathcal H$, by defining $T_n=\bigcap_{j\leq n}H_j$ for each $n\in\omega$, we obtain a decreasing sequence $(T_n)_{n\in\omega}$ in $\mathcal F\subseteq \mathcal I^+$.
By Proposition~\ref{prop:small-ideal}, there exists $I\in\mathcal I^+$ such that $I\subseteq^*T_n\subseteq H_n$ for all $n\in \omega$.
Now fix $\mathcal K\in[\mathcal H]^{<\omega}$.
As $\mathcal K$ is finite, $\bigcap\mathcal K\in\mathcal F\subseteq \mathcal I^+$, and therefore $I\not\subseteq^*(\omega\setminus\bigcap\mathcal K)\in \mathcal I$, which implies that $|I\setminus(\omega\setminus\bigcap\mathcal K)|=|I\cap\bigcap\mathcal K|=\omega$.
\end{proof}

\section{The construction}

\begin{definition}\label{def:maximal-centered}
Given an AD family $\mathcal A$ and a fin-sequence $C$,
we say that $\mathcal B\subseteq\mathcal A$ is \emph{maximal with
centered $C$-traces} if $\pi_C[\mathcal B]$ is centered and
$\pi_C[\mathcal B\cup\{a\}]$ is not centered for every
$a\in\mathcal A\setminus \mathcal B$.
\end{definition}

Equivalently, besides centeredness, we require that for
each $a\in\mathcal A\setminus \mathcal B$ there be
$F_a\in[\mathcal B]^{<\omega}$ such that
\begin{equation}\label{eq:finite-obstruction}
 |\pi_C(a)\cap J_C(F_a)|<\omega.
\end{equation}

We follow Corral and Rodrigues \cite[Theorem~5]{CorralRodrigues2024} by listing each pair $(C,X)$ cofinally often and updating the families
$B_C^\alpha$ by adding new members whenever their $C$-traces preserve
centeredness.

The preceding definition will be used to state Lemma~\ref{lem:local-extension}, a key result for our construction that allows us to carry out the recursion under $\mathfrak b=\mathfrak c$, which will be stated below and proved in the last section.

\begin{unnumberedtheorem}
    If $\mathfrak b=\mathfrak c$, then every infinite almost disjoint family of size less than $\mathfrak c$ can be extended to a pseudocompact MAD family.
\end{unnumberedtheorem}
\begin{proof}
Assume $\mathfrak b=\mathfrak c$ and fix an infinite AD family
$\mathcal D$ of size $\kappa<\mathfrak c$. Let $\mathscr C$
be the collection of all fin-sequences. Fix an enumeration
\[
 \bigl((C_\alpha,X_\alpha):\kappa\leq\alpha<\mathfrak c\bigr)
\]
of $\mathscr C\times\mathcal P(\omega)$ in which each pair
occurs cofinally often. Such an enumeration exists because
$|\mathscr C\times\mathcal P(\omega)|=\mathfrak c$ and
$\mathfrak c\cdot\mathfrak c=\mathfrak c$: arrange for
each pair to occur on a set of indices of size $\mathfrak c$.
Each such set is unbounded in the initial ordinal
$\mathfrak c$.

Enumerate $\mathcal D$ without repetitions as $(a_\xi)_{\xi<\kappa}$.
We recursively extend this enumeration to a sequence
$(a_\alpha)_{\alpha<\mathfrak c}$, simultaneously defining
\[
 \mathcal A_\alpha=\{a_\xi:\xi<\alpha\},
 \qquad
 B_C^\alpha\subseteq\mathcal A_\alpha
 \quad(\alpha\leq\mathfrak c,\ C\in\mathscr C).
\]
Start with $\mathcal A_0=B_C^0=\varnothing$ for every $C$.
Once $a_\alpha$ has been specified, set
$\mathcal A_{\alpha+1}=\mathcal A_\alpha\cup\{a_\alpha\}$
and, for every $C\in\mathscr C$, let
\begin{equation}\label{eq:greedy}
 B_C^{\alpha+1}=
 \begin{cases}
 B_C^\alpha\cup\{a_\alpha\},
   &\text{if }\pi_C[B_C^\alpha]\cup\{\pi_C(a_\alpha)\}
           \text{ is centered},\\
 B_C^\alpha,&\text{otherwise}.
 \end{cases}
\end{equation}
At a nonzero limit ordinal $\lambda\leq\mathfrak c$, put
\[
 \mathcal A_\lambda=\bigcup_{\xi<\lambda}\mathcal A_\xi,
 \qquad B_C^\lambda=\bigcup_{\xi<\lambda}B_C^\xi
 \quad(C\in\mathscr C).
\]
For $\alpha<\kappa$, use the prescribed $a_\alpha$.
These clauses determine $B_C^\kappa$ for every $C$.
At each later stage, the choice of $a_\alpha$ described
below depends only on objects already defined at stage
$\alpha$, after which \eqref{eq:greedy} is applied to all $C$.

We record the invariant furnished by this rule.
For each $\alpha\leq\mathfrak c$ and $C\in\mathscr C$,
$B_C^\alpha$ is maximal with centered $C$-traces in
$\mathcal A_\alpha$. Centeredness is preserved at successors
by definition and at limits because each finite subfamily
appears at an earlier stage. If $a_\xi$ was not included
in $B_C^{\xi+1}$, there was already a finite
$F\subseteq B_C^\xi$ such that
$|\pi_C(a_\xi)\cap J_C(F)|<\omega$.
That witness remains in every later candidate.
This proves the invariant by induction.

Now let $\kappa\leq\alpha<\mathfrak c$, and suppose the
construction has been carried out below $\alpha$.
Then $\mathcal A_\alpha$ is an infinite AD family of
cardinality less than $\mathfrak b$. Consider the following lemma:

\begin{lemma}\label{lem:local-extension}
Let $\mathcal A$ be AD with $|\mathcal A|<\mathfrak b$, $C$ a fin-sequence, and
$\mathcal B\subseteq\mathcal A$ maximal with centered $C$-traces.
Suppose that $X\subseteq\omega$ satisfies
$b\subseteq^*X$ for all $b\in \mathcal B$, and that for some
$F_0\in[\mathcal B]^{<\omega}$,
\[
 |J_C(F_0)\cap E_C(X)|<\omega.
\]
There exists $r\in[\omega]^\omega$ such that
\[
 r\perp\mathcal A,\qquad
 |r\setminus X|=\omega,\qquad
 \pi_C[\mathcal B]\cup\{\pi_C(r)\}\text{ is centered}.
\]
\end{lemma}

The hypotheses of the previous lemma motivate us to consider the following two conditions:

\begin{enumerate}[label=(\roman*), topsep=8pt]
 \item\label{item:upper-condition} $b\subseteq^*X_\alpha$
       for every $b\in B_{C_\alpha}^\alpha$;
 \item\label{item:bad-indices} there exists
       $F\in[B_{C_\alpha}^\alpha]^{<\omega}$ such that
       $|J_{C_\alpha}(F)\cap E_{C_\alpha}(X_\alpha)|<\omega$.
\end{enumerate}

If both hold, apply Lemma~\ref{lem:local-extension} to
$\mathcal A_\alpha,C_\alpha,B_{C_\alpha}^\alpha,X_\alpha$
and choose the resulting set as $a_\alpha$.
In particular,
\[
 a_\alpha\perp\mathcal A_\alpha,\qquad
 |a_\alpha\setminus X_\alpha|=\omega,\qquad
 a_\alpha\in B_{C_\alpha}^{\alpha+1}.
\]
If one of the conditions fails, choose any infinite
$a_\alpha\perp\mathcal A_\alpha$. Such a set exists because
$|\mathcal A_\alpha|<\mathfrak b\leq\mathfrak a$.
Update all candidates by \eqref{eq:greedy}.
This completes the recursion.

Set
\[
 \mathcal A=\mathcal A_{\mathfrak c}\qquad \text{and}
 \qquad B_C=B_C^{\mathfrak c}\quad(\text{for } C\in\mathscr C).
\]

The family $\mathcal A$ contains $\mathcal A_\kappa=\mathcal D$, is AD and has size $\mathfrak c$,
since every new infinite member is almost disjoint from
its predecessors and hence distinct from them.
For each $C$, the family $\pi_C[B_C]$ is centered.
Moreover, $B_C\neq\varnothing$: at a stage with
$(C_\alpha,X_\alpha)=(C,\varnothing)$, either
$B_C^\alpha$ is already nonempty, or
conditions \ref{item:upper-condition} and \ref{item:bad-indices} hold
(with $F=\varnothing$), and $a_\alpha$ is included.

We claim that $B_C$ is an accumulation point of $C$
for every $C\in\mathscr C$. If not, then by
Lemma~\ref{lem:vietoris}\ref{item:vietoris-criterion}, there exist
$X\subseteq\omega$ and finite $F\subseteq B_C$ such that
\[
 \forall b \in B_C\, (b\subseteq^*X)
 \qquad \text{and}\qquad |J_C(F)\cap E_C(X)|<\omega.
\]
Choose an occurrence $(C_\alpha,X_\alpha)=(C,X)$ after
all members of $F$ have been included in $B_C^\alpha$.
This is possible by the cofinal repetition of the pairs.
Then, both conditions \ref{item:upper-condition} and \ref{item:bad-indices} hold at stage $\alpha$. Thus $a_\alpha\in B_C$ and
$|a_\alpha\setminus X|=\omega$, a contradiction.

To check maximality, suppose $Y\in[\omega]^\omega$ were
almost disjoint from $\mathcal A$, and enumerate it
without repetitions as $(y_n)_{n\in\omega}$.
The fin-sequence $C_n=\{y_n\}$ has the nonempty accumulation
point $B_C\subseteq\mathcal A$. For any $a\in B_C$, its
lower Vietoris neighborhood determined by $\{a\}\cup a$
would imply $|Y\cap a|=\omega$, a contradiction.
Hence $\mathcal A$ is MAD.
Finally, Proposition~\ref{prop:finite-criterion} yields
the pseudocompactness of $\Exp(\Psi(\mathcal A))$.
\end{proof}

\section{Proof of the extension lemma}

It only remains to prove Lemma~\ref{lem:local-extension}. We shall do so by establishing some auxiliary lemmas.

\begin{lemma}\label{lem:selector}
Let $\mathcal A$ be AD with $|\mathcal A|<\mathfrak b$,
and let $C$ be a fin-sequence. Suppose that
$(a_i)_{i\in\omega}$ consists of distinct members of
$\mathcal A$ and that $\pi_C(a_i)$ is cofinite for every $i$.
There exist $m_0\in\omega$ and $w\perp\mathcal A$ such that
\[
 w\subseteq\bigcup_{m\geq m_0}C_m
 \quad\text{and}\quad
 |w\cap C_m|=1\quad(m\geq m_0).
\]
\end{lemma}
\begin{proof}
For each $i$, choose $r(i)$ such that, for every $m\geq r(i)$,
we can fix
\[
 x_{m,i}\in C_m\cap\left(a_i\setminus\bigcup_{j<i}a_j\right).
\]
This is possible as $a_i\cap\bigcup_{j<i}a_j$ is finite and $\pi_C(a_i)$ is cofinite.

For $a\in\mathcal A\setminus\{a_i:i\in\omega\}$, choose
$f_a\in\omega^\omega$ such that
$x_{m,i}\notin a$ whenever $m\geq\max\{r(i),f_a(i)\}$.
Such thresholds exist because $a\cap a_i$ is finite.

Choose a strictly increasing $q\in\omega^\omega$ which
eventually dominates all the $f_a$ and satisfies
$q\geq r$ pointwise. This is possible
since $|\mathcal A|<\mathfrak b$.
For $m\geq q(0)$, put
\[
 k(m)=\max\{i:q(i)\leq m\},\qquad
 w=\{x_{m,k(m)}:m\geq q(0)\}.
\]

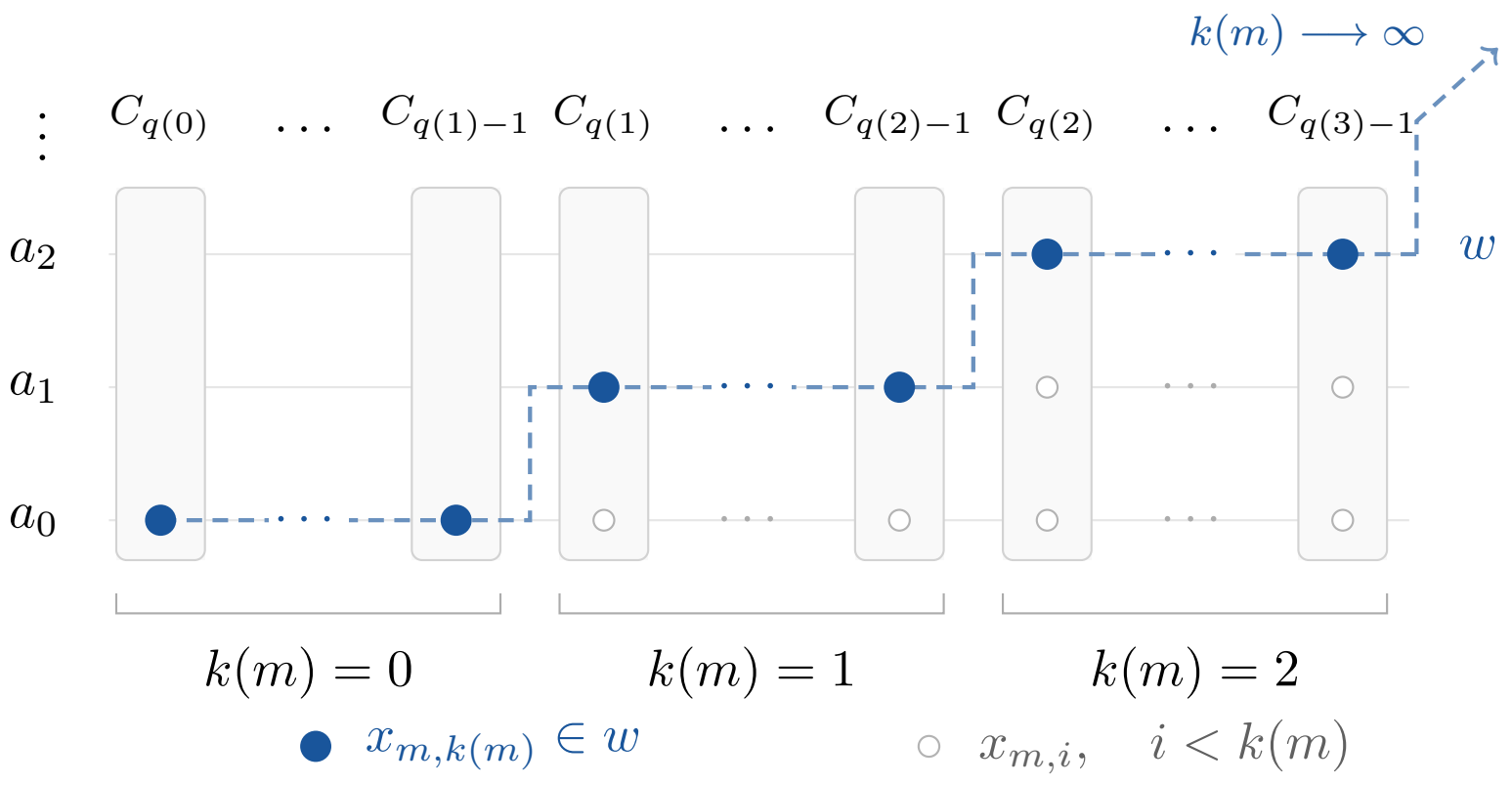
\begin{figure}[htbp]
\centering
\begin{tikzpicture}[x=1cm,y=0.9cm,font=\small]
  \definecolor{selectorblue}{RGB}{25,85,155}
  \foreach \y in {0,1,2}{
    \draw[gray!20] (-0.35,\y) -- (8.45,\y);
    \node[anchor=east] at (-0.55,\y) {$a_{\y}$};
  }
  \node at (-0.8,3) {$\vdots$};
  \foreach \x/\lab in {
    0/{C_{q(0)}},2/{C_{q(1)-1}},
    3/{C_{q(1)}},5/{C_{q(2)-1}},
    6/{C_{q(2)}},8/{C_{q(3)-1}}}{
    \fill[gray!5] (\x-0.3,-0.3) rectangle (\x+0.3,2.5);
    \draw[gray!35,rounded corners=2pt]
      (\x-0.3,-0.3) rectangle (\x+0.3,2.5);
    \node[font=\scriptsize,anchor=south] at (\x,2.7) {$\lab$};
  }
  \foreach \x in {1,4,7}{
    \node at (\x,2.93) {$\cdots$};
  }
  \foreach \x/\y in {3/0,5/0,6/0,8/0,6/1,8/1}{
    \draw[gray!60,fill=white] (\x,\y) circle (2pt);
  }
  \foreach \x/\y in {4/0,7/0,7/1}{
    \node[text=gray!65] at (\x,\y) {$\cdots$};
  }
  \draw[selectorblue!65,densely dashed,thick]
    (0,0) -- (2.5,0) -- (2.5,1) -- (5.5,1) -- (5.5,2) -- (8.5,2);
  \draw[selectorblue!65,densely dashed,thick,->]
    (8.5,2) -- (8.5,3) -- (9.05,3.55);
  \foreach \x/\y in {0/0,2/0,3/1,5/1,6/2,8/2}{
    \fill[selectorblue] (\x,\y) circle (3pt);
  }
  \foreach \x/\y in {1/0,4/1,7/2}{
    \node[fill=white,inner sep=1pt,text=selectorblue] at (\x,\y) {$\cdots$};
  }
  \node[anchor=west,text=selectorblue] at (8.65,2.05) {$w$};
  \node[anchor=east,text=selectorblue,font=\scriptsize]
    at (8.7,3.65) {$k(m)\longrightarrow\infty$};
  \foreach \left/\right/\middle/\i in {
    -0.3/2.3/1/0,2.7/5.3/4/1,5.7/8.3/7/2}{
    \draw[gray!65] (\left,-0.55) -- (\left,-0.7)
      -- (\right,-0.7) -- (\right,-0.55);
    \node[anchor=north] at (\middle,-0.8) {$k(m)=\i$};
  }
  \fill[selectorblue] (1.05,-1.7) circle (3pt);
  \node[anchor=west,text=selectorblue] at (1.25,-1.7) {$x_{m,k(m)}\in w$};
  \draw[gray!60,fill=white] (5.2,-1.7) circle (2pt);
  \node[anchor=west,text=gray!75!black] at (5.4,-1.7) {$x_{m,i},\quad i<k(m)$};
\end{tikzpicture}
\caption{$w=\{x_{m,k(m)}:m\geq q(0)\}$: one point per block.}
\label{fig:selector-points}
\end{figure}

Then $k(m)\to\infty$.
Given $j \in \omega$, notice that $x_{m, k(m)}\notin a_j$ whenever $k(m)>j$,
so $|w\cap a_j|<\omega$.
If $a$ is outside the chosen countable subfamily, then
$m\geq q(k(m))\geq f_a(k(m))$ for all sufficiently large $m$,
so the chosen point again avoids $a$.
Thus $w\perp\mathcal A$, and the conclusion follows with
$m_0=q(0)$.
\end{proof}

Although the next observation is a well-known fact, we shall include its proof for the convenience of the reader. It will be used along with the fact that $\mathfrak b\leq \mathfrak a$.

\begin{lemma}\label{lem:positive-orthogonal}
If $\mathcal A$ is AD, $|\mathcal A|<\mathfrak a$,
and $Y\in\mathcal I^+(\mathcal A)$, then some
$z\in[Y]^\omega$ satisfies $z\perp\mathcal A$.
\end{lemma}
\begin{proof}
Consider
$\mathcal A\restriction Y
=\{a\cap Y:a\in\mathcal A,\ |a\cap Y|=\omega\}$.
If this family is infinite, its size is less than
$\mathfrak a$, so it is not maximal on $Y$.
Any infinite subset of $Y$ almost disjoint from it
is almost disjoint from all of $\mathcal A$.
If it is finite, removing its union from $Y$ leaves
an infinite set almost disjoint from $\mathcal A$,
since $Y\in\mathcal I^+(\mathcal A)$.
\end{proof}
For the following lemma, recall that $\mathfrak b\leq \mathfrak d$.
\begin{lemma}\label{lem:centered-extension}
Let $\mathcal A$ be AD with $|\mathcal A|<\mathfrak b$,
let $C$ be a fin-sequence, and let $\mathcal B\subseteq\mathcal A$
be infinite. If $\pi_C[\mathcal B]$ is centered, then there exists
$w\in[\omega]^\omega$ such that $w\perp\mathcal A$ and
$\pi_C[\mathcal B]\cup\{\pi_C(w)\}$ is centered.
\end{lemma}
\begin{proof}
Choose distinct $b_i\in \mathcal B$, $i\in\omega$.
Apply Corollary~\ref{cor:positive-pseudo} to
$\mathcal H=\pi_C[\mathcal B]$ and $H_i=\pi_C(b_i)$, using
$|\mathcal H|<\mathfrak d$ to obtain $I\in[\omega]^\omega$ almost contained in every
$\pi_C(b_i)$ and meeting every finite intersection from
$\pi_C[\mathcal B]$ infinitely.

Reindex $(C_n)_{n\in I}$ in increasing order. Each $b_i$ meets all but
finitely many of these blocks.
Lemma~\ref{lem:selector}, applied to the entire family
$\mathcal A$, gives $w\perp\mathcal A$ with
$\pi_C(w)=^*I$. Hence $\pi_C[\mathcal B]\cup\{\pi_C(w)\}$ is centered.
\end{proof}


\begin{lemma}\label{lem:positive-residual}
Let $\mathcal A$ be AD, $C$ a fin-sequence, and
$\mathcal B\subseteq\mathcal A$ maximal with centered $C$-traces.
Suppose that $X\subseteq\omega$ satisfies
$b\subseteq^*X$ for all $b\in \mathcal B$, and that for some
$F_0\in[\mathcal B]^{<\omega}$,
\[
 |J_C(F_0)\cap E_C(X)|<\omega.
\]
Then, for every $F\in[\mathcal B]^{<\omega}$,
\[
 Y_F=\left(\bigcup_{n\in J_C(F)}C_n\right)\setminus X
 \in\mathcal I^+(\mathcal A).
\]
\end{lemma}
\begin{proof}
Suppose $Y_F\subseteq^*\bigcup K$ for some finite
$K\subseteq\mathcal A$. For $a\in K\setminus \mathcal B$, choose
$F_a$ as in \eqref{eq:finite-obstruction}, and put
\[
 G=F\cup F_0\cup\bigcup_{a\in K\setminus \mathcal B}F_a.
\]
The set $J_C(G)$ is infinite and has finite intersection
with $E_C(X)$. Since the blocks are disjoint, $Y_G$ is infinite.
For $a\in K\cap \mathcal B$, the intersection $Y_G\cap a$
is contained in the finite set $a\setminus X$.
For $a\in K\setminus \mathcal B$, only finitely many blocks indexed
by $J_C(G)$ meet $a$, by \eqref{eq:finite-obstruction}.
Thus $Y_G\cap a$ is finite for every $a\in K$, contradicting
$Y_G\subseteq Y_F\subseteq^*\bigcup K$.
\end{proof}

Finally, we restate and prove Lemma~\ref{lem:local-extension}.

\begin{unnumberedlemma}
Let $\mathcal A$ be AD with $|\mathcal A|<\mathfrak b$, $C$ a fin-sequence, and
$\mathcal B\subseteq\mathcal A$ maximal with centered $C$-traces.
Suppose that $X\subseteq\omega$ satisfies
$b\subseteq^*X$ for all $b\in \mathcal B$, and that for some
$F_0\in[\mathcal B]^{<\omega}$,
\[
 |J_C(F_0)\cap E_C(X)|<\omega.
\]
There exists $r\in[\omega]^\omega$ such that
\[
 r\perp\mathcal A,\qquad
 |r\setminus X|=\omega,\qquad
 \pi_C[\mathcal B]\cup\{\pi_C(r)\}\text{ is centered}.
\]
\end{unnumberedlemma}

\begin{proof}
If $\mathcal B$ is finite,
apply Lemma~\ref{lem:positive-orthogonal} to $Y_{\mathcal B}$.
The resulting $r\in[Y_{\mathcal B}]^\omega$ is almost disjoint from
$\mathcal A$, is disjoint from $X$, and satisfies
$\pi_C(r)\subseteq J_C(\mathcal B)$. Its trace is infinite because
the blocks are finite, so the required family is centered.

Suppose that $\mathcal B$ is infinite.
By Lemma~\ref{lem:centered-extension}, choose
$w\perp\mathcal A$ such that
$\pi_C[\mathcal B]\cup\{\pi_C(w)\}$ is centered.
Independently, Lemmas~\ref{lem:positive-residual}
and~\ref{lem:positive-orthogonal} give
$z\in[Y_\varnothing]^\omega$ with $z\perp\mathcal A$.
Let $r=w\cup z$. Then $r\perp\mathcal A$,
$z\subseteq r\setminus X$, and
$\pi_C(w)\subseteq\pi_C(r)$. The last inclusion preserves
all the required infinite intersections.
\end{proof}

\section*{Acknowledgements}

This study was financed, in part, by the S\~ao Paulo Research Foundation
(FAPESP), Brasil. Process Numbers 2025/07302-0 and 2026/06528-7.

\section*{Declaration of AI Assistance}
The main mathematical results presented in this paper had the first draft of their proofs generated by OpenAI's GPT-6 Astra. The authors carefully reviewed and verified the correctness of these proofs, rewriting them as necessary to ensure their accuracy and clarity. OpenAI's GPT-6.1 Sol, GPT-6.1 Terra, GPT-6.2 Sol and GPT-6.2 Astra aided the authors with subsequent revisions of the proofs, grammar, and LaTeX formatting.
\printbibliography

\end{document}